\documentclass[microtype]{gtpart}    
\usepackage{graphicx}  
\usepackage[all]{xy}
\usepackage{enumerate}
\usepackage{float}
\usepackage{mathrsfs} 
\usepackage{caption}
\title[The Chromatic Splitting Conjecture at $n=p=2$]{The Chromatic Splitting Conjecture at $n=p=2$}

\author[A. Beaudry]{Agn\`es Beaudry}
\givenname{Agn\`es}
\surname{Beaudry}
\address{Department of Mathematics\\ University of Colorado Boulder \\ \newline Campus Box 395 \\ Boulder \\ Colorado \\ 80309-0395}
\email{agnes.beaudry@colorado.edu}

\keyword{K(2)-local}
\keyword{chromatic}
\subject{primary}{msc2000}{55Q45}
\subject{primary}{msc2000}{55P60}

\arxivreference{math.AT/1412.2822}
\arxivpassword{6wbvm}

\volumenumber{}
\issuenumber{}
\publicationyear{}
\papernumber{}
\startpage{}
\endpage{}
\doi{}
\MR{}
\Zbl{}
\received{}
\revised{}
\accepted{}
\published{}
\publishedonline{}
\proposed{}
\seconded{}
\corresponding{}
\editor{}
\version{}

\newtheorem{thm}{Theorem}[subsection]
\newtheorem{cor}{Corollary}[subsection]
\newtheorem{prop}{Proposition}[subsection]
\newtheorem{lem}{Lemma}[subsection]
\theoremstyle{definition}
\newtheorem{defn}{Definition}[subsection]

\newtheorem{rem}{Remark}[subsection]
\newtheorem{notn}{Notation}[subsection]
\newtheorem{warn}{Warning}[subsection]

\makeatletter
\let\c@lem=\c@thm
\let\c@cor=\c@thm
\let\c@notn=\c@thm
\let\c@prop=\c@thm
\let\c@lem=\c@thm
\let\c@defn=\c@thm
\let\c@exmps=\c@thm
\let\c@rem=\c@thm
\let\c@warn=\c@thm
\makeatother

\numberwithin{equation}{subsection}

\newcommand{\F}{\mathbb{F}}
\newcommand{\Sn}{\mathbb{S}}
\newcommand{\G}{\mathbb{G}}
\newcommand{\W}{\mathbb{W}}
\newcommand{\smsh}{\wedge}
\newcommand{\ra}{\rightarrow}
\newcommand{\xra}{\xrightarrow}
\newcommand{\cC}{\mathcal{C}}

\DeclareFontFamily{OMS}{rsfs}{\skewchar\font'60}
\DeclareFontShape{OMS}{rsfs}{m}{n}{<-5>rsfs5 <5-7>rsfs7 <7->rsfs10 }{}
\DeclareSymbolFont{rsfs}{OMS}{rsfs}{m}{n}
\DeclareSymbolFontAlphabet{\scr}{rsfs}

\newcommand{\sC}{\mathscr{C}}
\newcommand{\sE}{\mathscr{E}}

\makeop{colim}
\makeop{Der}
\makeop{Gal}
\makeop{id}
\makeop{Mod}
\makeop{Tot}
\makeop{gr}
\makeop{Out}
\makeop{Hom}
\makeop{Ext}
\makeop{End}
\makeop{Aut}
\makeop{Tor}
\makeop{ev}
\makeop{maps}

\begin{document}

\begin{abstract}    
We show that the strongest form of Hopkins' chromatic splitting conjecture, as stated by Hovey in \cite{cschov}, cannot hold at chromatic level $n=2$ at the prime $p=2$. More precisely, for $V(0)$ the mod $2$ Moore spectrum, we prove that $\pi_{k}L_1L_{K(2)}V(0)$ is not zero when $k$ is congruent to $-3$ modulo $8$. We explain how this contradicts the decomposition of $L_1L_{K(2)}S$ predicted by the chromatic splitting conjecture.
\end{abstract}

\maketitle


\section{Introduction}

Fix a prime $p$. Let $S$ be the $p$-local sphere spectrum and $L_nS$ be the Bousfield localization of $S$ at the Johnson-Wilson spectrum $E(n)$. Let $K(n)$ be Morava $K$-theory. There is a homotopy pull-back square, called the chromatic fracture square, given by
 \begin{equation*}\label{hpp}
\xymatrix{L_{n}S \ar[r] \ar[d] & L_{K(n)}S \ar[d] \\ L_{n-1}S \ar[r]^-{\iota}&  L_{n-1}L_{K(n)}S.}\end{equation*}
Let $F_n$ be the fiber of the map $L_nS \ra L_{K(n)}S$. Note that $F_n$ is weakly equivalent to the fiber of $\iota$. It was shown in \cite[Lemma 4.1]{cschov} that $F_n$ is weakly equivalent to the function spectrum $F(L_{n-1}S, L_nS)$. Hopkins' chromatic splitting conjecture, as stated by Hovey \cite[Conjecture 4.2]{cschov}, stipulates that $\iota$
is the inclusion of a wedge summand, so that
\begin{equation}\label{splitting}L_{n-1}L_{K(n)}S \simeq L_{n-1}S \vee \Sigma F_n. \end{equation}
We will call this the \emph{weak} form of the chromatic splitting conjecture.
However, \cite[Conjecture 4.2]{cschov} also gives an explicit decomposition of $\Sigma F_n$ as a wedge of suspensions of spectra of the form $L_iS_p$ for $0\leq i<n$. We will call this the \emph{strong} form of the chromatic splitting conjecture.

The conjectured decomposition comes from the connection between the $K(n)$--local category and the cohomology of a certain group called the Morava stabilizer group $\G_n$. Let $\mathbb{S}_n$ be the group of automorphisms of the formal group law of $K(n)$ over $\F_{p^n}$. Then $\G_n$ is the extension of $\mathbb{S}_n$ by the Galois group $\Gal(\F_{p^n}/\F_p)$. Let $\W$ be the Witt vectors on $\F_{p^n}$. There is a spectral sequence
\begin{equation}\label{DSS}
H^s(\G_n, (E_n)_t) \Longrightarrow \pi_{t-s}L_{K(n)}S.
\end{equation}
Note that $\W$ sits naturally in $(E_n)_0 \cong \W[\![u_1, \ldots, u_{n-1}]\!]$. The inclusion induces a map
\begin{equation}\label{mapWEn}
H^*(\G_n, \W) \ra H^*(\G_n, (E_n)_0).\end{equation}
Morava proves in \cite[Remark 2.2.5]{morava} using the work of Lazard that
\[H^*(\G_n,\W)\otimes \Q_p \cong E(e_1, \ldots, e_n)\]
for classes $e_i$ of degree $2i-1$. Therefore, $H^*(\G_n,\W)$ contains an exterior algebra $E(x_1, \ldots, x_n)$ for appropriate integral multiples $x_i$ of the generators $e_i$. The chromatic splitting conjecture stipulates that, for some choice of $x_i$'s, the exterior algebra $E(x_1,\ldots, x_n)$ injects into $H^*(\G_n, (E_n)_0)$ under the map (\ref{mapWEn}), and that the non-zero classes $x_{i_1}\ldots x_{i_j}$ survive in (\ref{DSS}) to non-trivial elements in $\pi_{-2(\sum i_k)+j }L_{K(n)} S$. Further, it states that there is a factorization
 \[ \xymatrix{  S_p^{-2(\sum i_k)+j} \ar[r]  \ar[d] & L_{n-\max(i_k)}S_p^{-2(\sum i_k)+j} \ar@{.>}[d]  \\
  L_{K(n)}S \ar[r] & \Sigma F_n ,} \] 
where $S^m_p$ is the $p$--completion of $S^m$, and that these maps decompose $\Sigma F_n$ as
\begin{equation}\label{bigdecomp}
\Sigma F_n \simeq \bigvee_{\substack{1\leq j \leq n \\ 1 \leq i_1 < \ldots < i_j \leq n}} L_{n-\max(i_k)}S_p^{-2(\sum i_k)+j}.
\end{equation}

The chromatic splitting conjecture has been shown for $n\leq 2$ and for all primes $p$, except in the case $n=p=2$. For $n=1$, it follows immediately from a computation of $\pi_*L_1S_p$ (see Ravenel \cite[Theorem 8.10, 8.15]{ravloc}). At $n=2$ and $p\geq 5$, it is due to Hopkins, and follows from Shimomura's and Yabe's computations in \cite{shimyab}. The proof can be found in Behrens's account of their work \cite[Remark 7.8]{behse2}. At $n=2$ and $p=3$, the conjecture was proved recently by Goerss, Henn and Mahowald in \cite{GoerssSplit}. 

In this paper, we show that the chromatic splitting conjecture as stated above cannot hold for $n=p=2$. More precisely, we show that \cite[Conjecture 4.2(iv)]{cschov} fails in this case. At $n=2$, (\ref{splitting}) and (\ref{bigdecomp}) imply that
\begin{equation}\label{bigdecomp1}
L_{1}L_{K(2)}S \simeq L_1S_p \vee L_1S_p^{-1} \vee L_0 S_p^{-3}\vee L_0S_p^{-4}. 
\end{equation}
We show that the right hand side of (\ref{bigdecomp1}) has too few homotopy groups for the equivalence to hold. However, our results do not contradict the possibility that $\iota$ is the inclusion of a wedge summand. Giving an alternative description for the fiber in this case is work in progress. 

That our methods could disprove (\ref{bigdecomp1}) was first suggested to the author by Paul Goerss. He and Mark Mahowald had been studying the computations of Shimomura and Wang in \cite{shimwang} and \cite{shim} and noticed that these suggest that the right hand side of (\ref{bigdecomp1}) is too small.

\subsection*{Statement of the results}
\begin{figure}[H]
\center
\includegraphics[width=0.8\textwidth]{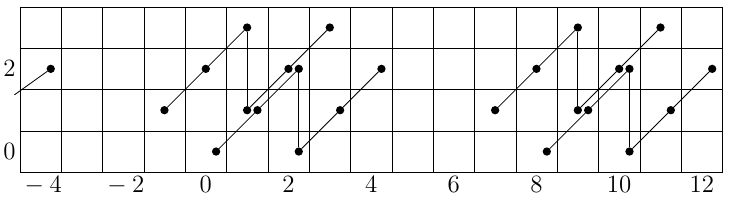}
\caption{The $E_{\infty}$--term of (\ref{DSS1V0}) computing $\pi_*L_{K(1)}V(0)$. Vertical lines denote extensions by multiplication by $2$ and lines of slope one denote multiplication by $\eta$.}
\label{fig:piL1}
\end{figure}

Let $V(0)$ be the cofiber of multiplication by $p$ on $S$. Note that for any $p$-local spectrum $X$, there is a cofiber sequence
\[X \xra{p} X \ra X \smsh V(0).\]
Since Bousfield localization of spectra preserves exact triangles, it follows that
\[L_{E}V(0) \simeq L_ES \smsh V(0)\]
for any spectrum $E$. This has the following consequence.
\begin{prop}\label{CSCV0}
The strong form of the chromatic splitting conjecture at $n=2$ implies that $L_{1}L_{K(2)}V(0) \simeq L_{1}V(0) \vee L_{1} \Sigma^{-1}V(0)$.
\end{prop}
We now fix our attention to the case when $p=2$. Since $L_{0}V(0)$ is contractible, it follows from the chromatic fracture square that 
$L_1V(0)\simeq L_{K(1)}V(0)$.
Computing $\pi_*L_{K(1)}V(0)$ is a routine exercise using the spectral sequence
\begin{equation}\label{DSS1V0}
E_{2}^{s,t} = H^s(\G_1,(E_1)_*V(0)) \Longrightarrow \pi_{t-s}L_{K(1)}V(0).
\end{equation}
The $E_{\infty}$--term is given in \fullref{fig:piL1}. At $p=2$, $V(0)$ is not a ring spectrum. This manifests itself by the fact that $\pi_*L_{K(1)}V(0)$ is not a ring. In fact,
\[\pi_*L_{K(1)}V(0) =  \left(\Z_2[\eta, \beta^{\pm 1},\zeta_1]/(2\eta, \eta^3, \zeta_1^2)\right)\{{e_0}, {v_1e_0}\}/(2  {e_0},2 {v_1}e_0 - \eta^2  {e_0}),\]
where $\eta \in \pi_1$ is the Hopf map, $\beta \in \pi_8$ is the $v_1$-self-map detected by $v_1^4$ and $\zeta_1 \in \pi_{-1}$ is detected by a generator of $H^1(\G_1; \Z_2) \cong H^1(\Z_2^{\times}, \Z_2)$. The element $e_0 \in \pi_0$ represents the inclusion of the bottom cell $S^0 \hookrightarrow V(0)$ and
$v_1e_0 \in \pi_2$ is a lift of $\Sigma \eta$ to the top cell:
\[\xymatrix{  & & S^2 \ar@{.>}[dl]_{v_1e_0} \ar[d]^{\Sigma \eta} & \\
S^0 \ar[r]^-{e_0} & V(0) \ar[r] & S^1 \ar[r]^2 &S^1. }\]
The following result is a consequence of \fullref{CSCV0}.
\begin{cor}
The chromatic splitting conjecture implies that $\pi_kL_1L_{K(2)}V(0)$ is zero when $k\equiv -3$ modulo $8$.
\end{cor}
However, in this paper, we prove the following result.

\begin{thm}\label{main}
There are non-trivial homotopy classes $\beta^t x $ in $\pi_{8t-3}L_1L_{K(2)}V(0) $ and $\zeta_2\beta^t x $ in $\pi_{8t-4}L_1L_{K(2)}V(0)$.
\end{thm}
This has the following immediate consequence.
\begin{thm}
The homotopy group $\pi_{k}L_1L_{K(2)}V(0)$ is non-zero when $k \equiv -3$ modulo $8$. Therefore, the decomposition (\ref{bigdecomp1}) of the chromatic splitting conjecture does not hold when $n=2$ and $p=2$.
\end{thm}

The broad strokes of the proof of \fullref{main} when $t=0$ are as follows. Let $G_{24} \cong Q_8 \rtimes C_3$ be a representative of the unique conjugacy class of maximal finite subgroups of $\mathbb{S}_2$. Let $C_6 $ be a subgroup of $G_{24}$ of order $6$. Let $\mathbb{S}_2^1$ be the norm one subgroup so that $\mathbb{S}_2 \cong \mathbb{S}_2^1 \rtimes \Z_2$ (see \fullref{sec:e1loc}). It follows from the duality resolution techniques of Goerss, Henn, Mahowald and Rezk and the work of Bobkova in \cite{Bobkova} that, for any $X$, there is a spectral sequence
\begin{equation*}
\xymatrix{E_1^{p,t} = \pi_{t}(\sE_p\smsh X)  \ar@{=>}[r] & \pi_{t-p}(E_{2}^{h\mathbb{S}_{2}^1}\smsh X),} 
\end{equation*}
where $\sE_p$ are spectra such that $\sE_0 \simeq  E_2^{hG_{24}} $, $\sE_p  \simeq  E_2^{hC_{6}}$ if $p=1,2$ and $(E_{\cC})_{*}\sE_3  \cong   (E_{\cC})_*E_2^{hG_{24}}$ as Morava modules.
Localizing at $E(1)$, we obtain a spectral sequence
\begin{equation}\label{SSintro}
\xymatrix{E_1^{p,t} = \pi_{t}L_1(\sE_p\smsh X)  \ar@{=>}[r] & \pi_{t-p}L_1(E_{2}^{h\mathbb{S}_{2}^1}\smsh X).} 
\end{equation}
In \fullref{lem:e23} and \fullref{pi5S21}, we use this spectral sequence to show that
$\pi_{-3}L_1(E_{2}^{h\mathbb{S}_{2}^1} \smsh V(0)) \cong \F_4$. After taking Galois invariants, we obtain 
a non-zero element $x$ in $\pi_{-3}L_1(E_{2}^{h\G_{2}^1} \smsh V(0))$. In the cofiber sequence
\[L_1L_{K(2)}V(0)  \ra L_1(E_{2}^{h\G_{2}^1} \smsh V(0)) \ra L_1(E_{2}^{h\G_{2}^1} \smsh V(0)), \]
which is obtained from the cofiber sequence $L_{K(2)}S \to E_{2}^{h\G_{2}^1}  \to E_{2}^{h\G_{2}^1} $ by smashing with $V(0)$ and localizing at $E(1)$, this class gives rise to non-zero elements $x\in  \pi_{-3}L_1L_{K(2)}V(0) $ and $\zeta_2x \in \pi_{-4} L_1L_{K(2)}V(0)$. 
\begin{warn}
We use the notation $\zeta_2$ to denote the homotopy class defined by 
\[\xymatrix{  & & S^0 \ar[d]_-{1} \ar[dr]^-{\zeta_2} & \\
L_{K(2)}S \ar[r] & E_{2}^{h\G_{2}^1}  \ar[r] & E_{2}^{h\G_{2}^1}  \ar[r] &  \Sigma L_{K(2)}S }\] 
Experts will notice that this clashes with Ravenel \cite[Lemma 2.1]{RavCoh}, but this is the natural generalization of what is now commonly denoted by $\zeta_n$ at odd primes.
\end{warn}

\subsection*{Organization of the paper}

In Section \ref{sec:e1loc}, we specialize to the case $n=2$ and $p=2$ and describe the duality resolution spectral sequence and its $E(1)$-localization. In Section \ref{sec:e1page}, we compute the $E_1$-page of this spectral sequence for $V(0)$. In Section \ref{sec:main}, we prove \fullref{main}.

\subsection*{Acknowledgements} 
I thank Paul Goerss, Hans-Werner Henn and Peter May for their constant help and support. I thank Tobias Barthel, Daniel G. Davis and Douglas Ravenel for helping me sort out some of the details for the proofs of \fullref{sec:e1loc}. I also thank Mark Behrens, Irina Bobkova, Michael Hopkins, Jack Morava, Niko Naumann and Zhouli Xu for helpful conversations. Finally, I thank Mark Mahowald for the insight he shared with all of us throughout his life.

This material is based upon work supported by the National Science Foundation under Grant No. DMS-1612020.

\section{The $E(1)$-local Duality Resolution Spectral Sequence}\label{sec:e1loc}

We take the point of view that, at height $2$, the Honda formal group law may be replaced by the formal group law of a supersingular elliptic curve. This was carefully explained in \cite[Section 1]{Paper2}. (The reader who wants to ignore this subtlety may take $\mathbb{S}_{\cC}$, $\G_{\cC}$ and $E_{\cC}$ to mean $\mathbb{S}_2$, $\G_2$ and $E_2$ respectively.)  

Let $\mathbb{S}_{\cC}$ be the group of automorphisms of the formal group law of the supersingular elliptic curve 
\[\cC: y^2+y=x^3\] 
of height two over $\F_4$ (see \cite[Section 3]{Paper2} for the comparison). It admits an action of the Galois group $\Gal(\F_4/\F_2)$. Define 
\[\G_{\cC} = \mathbb{S}_{\cC} \rtimes \Gal(\F_4/\F_2).\]
Let $E_{\cC}$ be the spectrum which classifies the deformations of the formal group law of $\cC$ over $\F_4$ as described, for example, in Rezk \cite{rezk}. It can be chosen to be a complex oriented ring spectrum with 
\[(E_{\cC})_* = \W[\![u_1]\!][u^{\pm1}]\]
for $|u_1|=0$, $|u|=-2$, whose formal group law is the formal group law of the curve
\begin{equation}\label{CU}
\cC_{U} : y^2 + 3u_1xy+(u_1^3-1)y = x^3.
\end{equation}
It admits an action of $\G_{\cC}$ and for any finite spectrum $X$,
\[L_{K(2)}X \simeq E_{\cC}^{h\G_{\cC}}\smsh X \simeq (E_{\cC} \smsh X)^{h\G_{\cC}} \]
(see Behrens and Davis \cite[p. 5]{behdav}).
The group of automorphisms $\Aut(\cC)$ of $\cC$ is of order $24$ and injects into $\mathbb{S}_{\cC}$. We let $G_{24}$ denote the image of $\Aut(\cC)$. We note that 
\[G_{24} \cong Q_8 \rtimes C_3,\]
where $Q_8$ is a quaternion subgroup and $C_3$ a cyclic group of order $3$. The group $\mathbb{S}_{\cC}$ contains a central subgroup of order $2$, which we denote by $C_2$. We define
\[C_6 = C_2 \times C_3.\]

There is a surjective homomorphism $N: \mathbb{S}_{\cC} \to  \Z_2^{\times}/(\pm 1) \cong \Z_2$, which we call the \emph{norm}. It is constructed using the determinant of a representation $\rho: \mathbb{S}_{\cC}  \to GL_2(\W)$ (see \cite[Section 3]{Paper2}). Further, it can be extended to $\G_{\cC}$. We let $\mathbb{S}_{\cC}^1$ and $\G_{\cC}^1$ be the kernels of the norms, and note that the elements of finite order in $\mathbb{S}_{\cC}$ and $\G_{\cC}$ are contained in $\mathbb{S}_{\cC}^1$ and $\G_{\cC}^1$ respectively. Further,
\begin{align}\label{semi}
\mathbb{S}_{\cC} &\cong \mathbb{S}_{\cC}^1 \rtimes \Z_2,  & \G_{\cC} &\cong \G_{\cC}^1 \rtimes \Z_2.  \end{align}

The formal group law $F_{\cC_U}$ of $\cC_U$, is not $2$--typical. Nonetheless, it is strictly isomorphic to a $2$--typical formal group law classified by a map $BP_* \to (E_{\cC})_*$. Further, $[2]_{F_{\cC_U}}(x) \equiv u_1 u^{-1}x^2$ modulo $(2,x^4)$ (see \cite[Section 6.1]{Paper2} for details on $F_{\cC_U}$). The strict isomorphism between $F_{\cC_U}$ and its $2$--typification preserves this identity. Hence, $v_1$ is mapped to $u_1u^{-1}$ modulo $(2)$. Since we are working primarily modulo $(2)$, we abuse notation and let
$v_1 = u_1 u^{-1} \in (E_{\cC})_2$.

We will need the following result, which can be found in Henn \cite[Theorem 13]{hennres} and is also discussed in greater details in Bobkova \cite{Bobkova}. We restate it here using our notation for convenience. 
\begin{thm}[Goerss, Henn, Mahowald, Rezk, Bobkova]\label{bob}
There is a resolution of spectra in the $K(2)$-local category given by
\[\xymatrix @R=1.3pc {E_{\cC}^{h\mathbb{S}_{\cC}^1} \ar[r] &  E_{\cC}^{hG_{24}} \ar[r] \ar@{=}[d] &  E_{\cC}^{hC_6} \ar[r]  \ar@{=}[d] & E_{\cC}^{hC_6} \ar[r]  \ar@{=}[d] & \sE_3 \ar@{=}[d] \\  
& \sE_0 \ar[r] & \sE_1 \ar[r] & \sE_2 \ar[r] & \sE_3} \]
where $(E_{\cC})_*\sE_3 \cong (E_{\cC})_*E_{\cC}^{hG_{24}}$ as Morava modules.
Further, for any spectrum $X$, the resolution gives rise to a tower of fibrations spectral sequence 
\begin{equation}\label{SS1}
\xymatrix{E_1^{p,t} = \pi_{t}(\sE_p\smsh X)  \ar@{=>}[r]^-{SS_1} & \pi_{t-p}(E_{\cC}^{h\mathbb{S}_{\cC}^1}\smsh X)} 
\end{equation}
with differentials $d_r : E_r^{p,t}  \ra E_r^{p+r,t+r-1}$.
\end{thm}

We call the resolution of \fullref{bob} the \emph{duality resolution}. Let $\pi$ generate $\Z_2$ in the decompositions (\ref{semi}) and let $G_{24}' = \pi G_{24} \pi^{-1}$. Recall from \cite{Paper2} or \cite{Paper1} that there is also an \emph{algebraic} duality resolution:
\begin{align}\label{adr}
\xymatrix @C=1.3pc @R=1.3pc { 0 \ar[r] & \Z_2[\![\mathbb{S}_{\cC}^1/G_{24}']\!]  \ar[r]  \ar@{=}[d] & \Z_2[\![\mathbb{S}_{\cC}^1/C_6]\!]  \ar[r]  \ar@{=}[d]  & \Z_2[\![\mathbb{S}_{\cC}^1/C_6]\!]   \ar[r]   \ar@{=}[d]  & \Z_2[\![\mathbb{S}_{\cC}^1/G_{24}]\!]  \ar[r]   \ar@{=}[d] & \Z_2 \ar[r] & 0 \\
&  \sC_3  \ar[r]  & \sC_2 \ar[r] & \sC_1  \ar[r] & \sC_0 & & }
\end{align}

Now, let $X$ be a finite spectrum. Resolving (\ref{adr}) into a double complex of projective $\mathbb{S}_{\cC}^1$--modules and applying the functor $\Hom_{\Z_2[\![\mathbb{S}_{\cC}^1]\!]}(-, ( E_{\cC})_tX)$ gives rise to a spectral sequence
 \begin{equation}\label{SS2}
 \xymatrix{E_1^{p,q,t} = \Ext^q_{\Z_2[\![\mathbb{S}_{\cC}^1]\!]}(\sC_p,( E_{\cC})_tX) \ar@{=>}[r]^-{SS_2} & H^{p+q}(\mathbb{S}_{\cC}^1,(E_{\cC})_tX)}\end{equation}
 with differentials $d_r : E_r^{p,q,t} \ra E_r^{p+r,q-r+1,t}$.
 Further, in each fixed degree $p$, there are spectral sequences
  \begin{equation}\label{SS3}
E_1^{p,q,t} = \xymatrix{ \Ext^q_{\Z_2[\![\mathbb{S}_{\cC}^1]\!]}(\sC_p, (E_{\cC})_tX)  \ar@{=>}[r]^-{SS_3} & \pi_{t-q}(\sE_p\smsh X)}  \end{equation}
 with differentials $d_r : E_r^{p,q,t} \ra E_r^{p,q+r,t+r-1}$.
 Finally, there is also a spectral sequence
    \begin{equation} \label{SS4}
 \xymatrix{ 
E_2^{s,t} =  H^{s}(\mathbb{S}_{\cC}^1,(E_{\cC})_tX)   \ar@{=>}[r]^-{SS_4} & \pi_{t-s}(E_{\cC}^{h\mathbb{S}_{\cC}^1}\smsh X)  }  \end{equation}
 with differentials $d_r : E_r^{s,t} \ra E_r^{s+r,t+r-1}$. Thus, for $X$ finite, we obtain a diagram of spectral sequences
  \begin{equation}\label{diagss} \xymatrix{ \Ext^q_{\Z_2[\![\mathbb{S}_{\cC}^1]\!]}(\sC_p, (E_{\cC})_tX)\ar@{=>}[d]_-{SS_3} \ar@{=>}[r]^-{SS_2} & H^{p+q}(\mathbb{S}_{\cC}^1,(E_{\cC})_tX) \ar@{=>}[d]^-{SS_4} \\
  \pi_{t-q}(\sE_p\smsh X)   \ar@{=>}[r]_-{SS_1} & \pi_{t-(p+q)}(E_{\cC}^{h\mathbb{S}_{\cC}^1}\smsh X). } \end{equation}

\begin{rem}\label{rem:edgehom}
For elements of Adams-Novikov filtration $s=0$ in $E_{1}^{p,t}(SS_1)$, the differentials $d_1$ are related to the $d_1$--differentials in the algebraic duality resolution spectral sequence $SS_2$ in the following way. If $X$ is finite, as in \cite[Proposition 2.4 and (2.7)]{ghmr}, for $G$ a closed subgroup of $\G_{\cC}$, there are isomorphisms of Morava modules
\begin{align}\label{eqn:ECmod}
 (E_{\cC})_t(E_{\cC}^{hG}\smsh X) \cong  \Hom^c(\G_{\cC}/G, (E_{\cC})_tX)  \cong  \Hom_{\Z_2}(\Z_2[\![\G_{\cC}/G]\!] , (E_{\cC})_tX)  .\end{align}
Let
\begin{equation*}
 \xymatrix{E_1(SS_1)^{p,t} \cong \pi_t(\sE_p \smsh X) \ar[r]^-{h} & H^0(\G_{\cC}, (E_{\cC})_t(\sE_p \smsh X)) \cong E_1^{p,0,t}(SS_2)}\end{equation*}
be the edge homomorphism for the spectral sequence 
\[H^s(\G_{\cC}, (E_{\cC})_t(\sE_p \smsh X)) \Longrightarrow \pi_{t-s}(\sE_p \smsh X). \]
The spectral sequence $SS_1$ is constructed so that the following diagram commutes:
\begin{equation*}
\xymatrix{ E_1^{p,t}(SS_1) \ar[r]^-{h} \ar[d]_-{d_1} & E_1^{p,0,t}(SS_2) \ar[d]^-{d_1} \\
  E_1^{p+1,t}(SS_1) \ar[r]^-{h} & E_1^{p+1,0,t}(SS_2)}.
  \end{equation*}
 When both horizontal maps $h$ are injective, one can deduce information in $SS_1$ from information in $SS_2$.
 \end{rem}
 
 For the statement of the next result, recall that for any closed subgroup $F$ of $\G_{\cC}$ and finite spectrum $X$ there is a spectral sequence
 \begin{equation}\label{eqn:descFG}
 E_{2}^{s,t}(F, X) = H^s(F, (E_{\cC})_tX ) \Longrightarrow \pi_{t-s} (E_{\cC}^{hF}\smsh X).\end{equation}
 The author has learned the proof of the following result from Paul Goerss.
 \begin{lem}\label{lem:Wlin}
Let $S$ a closed subgroup of $\mathbb{S}_{\cC}$ which is invariant under the action of $\Gal(\F_4/\F_2)$. Let $G \cong  S \rtimes \Gal(\F_4/\F_2)$ be the corresponding closed subgroup of $\G_{\cC}$. Then, for any finite $X$ and any $2\leq r\leq \infty$,
\[E_{r}^{s,t}(S, X) \cong \W \otimes_{\Z_2} E_{r}^{s,t}(G, X),  \]
and the differentials of the spectral sequence $E_{r}^{s,t}(S, X)$ are $\W$--linear.
\end{lem}
 \begin{proof}
The action of $\Gal(\F_4/\F_2)$ on $(E_{\cC})_*X$ is semilinear over $\W$, so there is an isomorphism $E_2^{*,*}(S,X)\cong  \W \otimes_{\Z_2}E_2^{*,*}(G,X)$. Now consider, $E_{r}^{s,t}(\Sn_{\cC}, S^0)$. We have $E_2^{0,0}(\Sn_{\cC}, S^0) \cong \W$ and the subring $\Z_2$ of $\W$ consists of permanent cycles. The spectral sequence $E_r^{*,*}(\Sn_{\cC},S^0)$ is multiplicative, so the differentials $d_r:  E_r^{0,0} \to E_r^{r,r-1}$ are $\Z_2$--derivations. Since $\W$ is an \'etale extension of $\Z_2$,  for any $r$, the $\Z_2$--derivations from $\W$ to the $\W$--module $E_r^{r,r-1}$ are zero.
Hence, $E_2^{0,0}(\Sn_{\cC}, S^0) \cong \W$ consists of permanent cycles and the differentials are $\W$--linear. 
Since the spectral sequence $E_r^{*,*}(S,X)$ is one of modules over $E_r^{*,*}(\Sn_{\cC},S^0)$, the differentials of $E_r^{*,*}(S,X)$ are also $\W$--linear and the result follows.
\end{proof}

In what follows, we will use the following remark.
\begin{rem}\label{rem:van}
Let $X$ be a finite spectrum and $F$ be a closed subgroup of $\G_{\cC}$.
As noted by Devinatz in the proof of \cite[Lemma 3.5]{devlshss}, it follows from the fact that $E_{\cC}^{hF}$ is $(K_{\mathcal{C}})_*$--local $E_{\mathcal{C}}$--nilpotent, (see Devinatz and Hopkins \cite[Proposition A.3]{DH}) that the descent spectral sequence (\ref{eqn:descFG})
has a horizontal vanishing line.
 \end{rem}

Now, recall that the Telescope conjecture holds at height $n=1$. This was proved at odd primes by Miller \cite{miller} and at $p=2$ by Mahowald \cite{mahimj}. In particular, we have the following result. 
 \begin{thm}[Mahowald, Miller]\label{telescope}
Let $Y$ admit a $v_1$ self-map $v_1^k : \Sigma^{2k} Y \ra Y$. Then
\begin{align*}
L_1Y \simeq L_{K(1)}Y \simeq v_1^{-1}Y,
\end{align*}
where
\begin{align*}
v_1^{-1}Y := \colim \left(\ldots \xrightarrow{v_1^k} \Sigma^{2k}Y \xrightarrow{v_1^k} Y \xrightarrow{v_1^k} \ldots \right).
\end{align*}
\end{thm}

 \begin{prop}\label{locss}
 For any finite type $1$ spectrum $X$, with self map $v_1^k : \Sigma^{2k}X \xrightarrow{v_1^k} X $, there is a diagram of strongly convergent spectral sequences
  \begin{equation*}\label{diagssv1} \xymatrix{ v_1^{-1}\Ext^q_{\Z_2[\![\mathbb{S}_{\cC}^1]\!]}(\sC_p, (E_{\cC})_tX)\ar@{=>}[d]_-{L_1SS_3} \ar@{=>}[rr]^-{L_1SS_2} & & v_1^{-1}H^{p+q}(\mathbb{S}_{\cC}^1,(E_{\cC})_tX) \ar@{=>}[d]^-{L_1SS_4} \\
  \pi_{t-q}L_1(\sE_p\smsh X)   \ar@{=>}[rr]_-{L_1SS_1} & & \pi_{t-(p+q)}L_1(E_{\cC}^{h\mathbb{S}_{\cC}^1}\smsh X). } \end{equation*}
 \end{prop}
\begin{proof}
The spectral sequence $L_1SS_2$ is obtained from $SS_2$ by inverting the element $v_1^k \in (E_{\cC})_{2k}X$, and $L_1SS_1$ is obtained by the applying $L_1$ to the tower of fibrations which gives rise to $SS_1$. 
The spectral sequences $L_1SS_3$ and $L_1SS_4$ are obtained by inverting the algebraic element $v_1^k$ in the spectral sequences $SS_3$ or $SS_4$, and using the fact that
\[v_1^{-1}\pi_{*}(\sE_p\smsh X)   \cong \pi_{*}L_1(\sE_p\smsh X). \]

With regards to the strong convergence of the four spectral sequences, note that localization with respect to $v_1$ is exact. Therefore, the localized spectral sequences will converge strongly if they have horizontal vanishing lines at the $E_{\infty}$--term. The spectral sequences $SS_1$ and $SS_2$ have a vanishing line at $p=4$ for all $r\geq 1$. 
As noted in \fullref{rem:van}, the descent spectral sequences $SS_3$ and $SS_4$ have horizontal vanishing lines. Therefore, the spectral sequences $L_1SS_i$ exist and converge.
 \end{proof}
 
 \begin{rem}\label{rem:edgehoml1}
 As in \fullref{rem:edgehom}, the differentials $d_1$ in $L_1SS_1$ and $L_1SS_2$ commute with the edge homomorphisms 
 \begin{equation*}
 \xymatrix{E_1(L_1SS_1)^{p,t} \cong \pi_tL_1(\sE_p \smsh X) \ar[r]^-{h} & v_1^{-1}H^0(\G_{\cC}, (E_{\cC})_t(\sE_p \smsh X)) \cong E_1^{p,0,t}(L_1SS_2)}\end{equation*}
 \end{rem}
 \begin{rem}\label{rem:Wlinl1}
For $X$ as in \fullref{locss}, the element $v_1^{2k} \in (E_{\cC})_{2k}X$ can be chosen to be Galois invariant. Therefore, the results of \fullref{lem:Wlin} also hold for the localized spectral sequences. That is, let
\[v_1^{-1}E_{2}^{s,t}(F, X) \cong v_1^{-1}H^s(F, (E_{\cC})_tX) \Longrightarrow \pi_{t-s}L_1(E_{\cC}^{hF}\smsh X).\]
Then for $S$ and $G$ are as in \fullref{lem:Wlin}, we have 
 \[ v_1^{-1}E_{r}^{s,t}(S, X) \cong \W \otimes_{\Z_2}v_1^{-1}E_{r}^{s,t}(G, X) \]
for $2\leq r \leq \infty$ and the differentials are $\W$--linear.
 \end{rem}

\section{The homotopy of $L_1(E_{\cC}^{hG_{24}}\smsh V(0))$ and $L_1(E_{\cC}^{hC_6}\smsh V(0))$}\label{sec:e1page}
The spectrum $V(0)$ has a self map
\[\beta  : \Sigma^{8}V(0)\xra{v_1^4} V(0),\]
and, in this section, we give the $E_1$--term for
\[\xymatrix{E_1^{p,q}(L_1SS_1) = \pi_{q}L_1(\sE_p\smsh V(0))  \ar@{=>}[r]^-{L_1SS_1} & \pi_{q-p}L_1(E_{\cC}^{h\mathbb{S}_{\cC}^1}\smsh V(0))}. \]
In order to do so, we must compute $\pi_*L_1(E_{\cC}^{hG_{24}}\smsh V(0))$ and $\pi_*L_1(E_{\cC}^{hC_6}\smsh V(0))$. We do this using the descent spectral sequences
\[v_1^{-1}H^s(G, (E_{\cC})_tV(0)) \Longrightarrow \pi_{t-s} L_1(E_{\cC}^{hG}\smsh V(0)).\] 

\begin{notn}
We use the following conventions. First,
\begin{align*}
v_1 &= u_1u^{-1},  & v_2 &= u^{-3}, & j_0 &= u_1^3.
\end{align*}
The element $\Delta$ is the discriminant of $\cC_{U}$, and hence is given by
\[\Delta = 27 v_2(v_1^3-v_2)^3 \equiv v_2(v_1^3+v_2)^3 \mod (2),\]
and
\[c_4 =  9 v_1^4+72 v_1 v_2 \equiv v_1^4 \mod (2).\]
The $j$-invariant is 
\[ j = c_4^3\Delta^{-1} \equiv v_1^{12}\Delta^{-1} \mod (2).\]
These identities can be computed using Silverman \cite[Section III.1]{silverman} (see also \cite[Section 4.2]{Paper2}).
We abuse notation and let
\[\eta   = \delta(v_1),\] 
where $\delta$ is the Bockstein associated to 
\[  0 \ra (E_{\cC})_*/2 \xra{2} (E_{\cC})_*/4 \ra (E_{\cC})_*/2  \ra 0. \]
This is justified by the fact that $\delta(v_1)$ detects the homotopy class $\eta$ (see \cite[Section 4.1]{Paper2}).
\end{notn}

The $v_1$--torsion free elements of $H^*(G_{24}, (E_{\cC})_*V(0))$ generate a submodule isomorphic to
\[  \F_4[\![j]\!][v_1,\eta, \Delta^{\pm 1}, k]/( \eta^4-v_1^4 k, j\Delta - v_1^{12})\]
for elements of degrees $(s,t)$, for $s$ the cohomological grading and $t$ the internal grading,
\begin{align*}
|v_1|&=(0,2), &  |\eta|  &= (1,2), & |\Delta|&=(0, 24), & |k|&=(4,0), & |j|&=(0,0),
\end{align*}
(see Section 4.2 or the Appendix of \cite{Paper2}).
On the other hand, $H^*(C_{6}, (E_{\cC})_*V(0))$ is $v_1$--torsion free and is isomorphic to
\begin{align*}
\mathbb{F}_4[\![j_0]\!][v_1,\eta, v_2^{\pm 1}, h]/(\eta-v_1 h,  j_0v_2 - v_1^3 ),
\end{align*}
where $|v_2| =(0,6)$, $|h| =(1,0)$ and $|j_0| =(0,0)$ (see Section 4.2 of \cite{Paper2}).

The next proposition is an immediate consequence of these results. In its statement, we let $\F_4(\!(x)\!)$ denote the Laurent series on $x$.
\begin{prop}\label{TBG24}
There are isomorphisms
\begin{align*}
v_1^{-1}H^*(G_{24}, (E_{\cC})_*V(0)) &\cong \F_{4}(\!(j)\!) [v_1^{\pm 1}, \eta ]
\end{align*}
and
\begin{align*}
v_1^{-1}H^*(C_{6} ; (E_{\cC})_*V(0)) &\cong  \F_{4}(\!(j_0)\!) [v_1^{\pm 1}, \eta ].
\end{align*}
The degrees $(s,t)$ are given by $|v_1 | =(0,2)$, $|\eta | =(1,2)$, $|j|=(0,0)$ and $|j_0| =(0,0)$. The restriction associated to the inclusion of $C_6$ in $G_{24}$ maps
$j$ to $j_0^{4}(1+j_0)^{-3}$.
\end{prop}
\begin{proof}
This follows from \cite[Section 4.2]{Paper2} after inverting $v_1$.
\end{proof}

To compute the differentials, we will use the following observation.
\begin{rem}\label{rem:unid3}
There is a class $\alpha_3$ in $\Ext_{BP_*BP}^{1,6}(BP_*,BP_*)$ (see Ravenel \cite[p.430]{ravnovice}) such that
$d_3(\alpha_3) = \eta^4$. Further, $\alpha_3$ reduces to $\eta v_1^2$ in $ \Ext_{BP_*BP}^{1,6}(BP_*,BP_*V(0))$, so that $\eta d_3(v_1^2) = \eta^4$. 

In general, for a $2$--local $BP$--algebra spectrum $E$, the $E$--Adams spectral sequence for any spectrum $X$ is a module over $\Ext_{BP_*BP}(BP_*,BP_*)$. There is a universal $d_3$--differential $d_3(\alpha_3 z) = \eta^4 z + \alpha_3 d_3(z)$.
Further, if $2$ annihilates $E_*(X)$, this reduces to $d_3(\eta v_1^2 z) = \eta^4 z +  \eta v_1^2 d_3(z)$. If there is no $\eta$--torsion on the $E_3$--term as in our examples below, this gives a universal differential
$d_3(v_1^2 z) = \eta^3 z +  v_1^2 d_3(z)$.
\end{rem}

\begin{lem}\label{lem:sse2}
Let $G$ be a closed subgroup of $\G_{\cC}$. Let $X$ be a $K(2)$--local spectrum such that $(E_{\cC})_*X \cong (E_{\cC})_*E_{\cC}^{hG}$. Then, the $K(2)$--local, $E_{\cC}$--Adams spectral sequence computing $\pi_{*}X$ has $E_2$--term isomorphic to $H^*(G, (E_{\cC})_*)$.
\end{lem}
\begin{proof}
We first prove that the $E_2$--term is isomorphic to $H^*(\G_{\cC}, (E_{\cC})_*X)$. This can be deduced directly from Barthel and Heard \cite[Theorem 4.3]{bartheldrew}. 
Nonetheless, we sketch the proof here. The assumption on $(E_{\cC})_*X$ implies that it is pro-free as an $(E_{\cC})_*$--module. An inductive argument using \cite[Proposition 8.4]{666} and \cite[Proposition 2.4]{ghmr} shows that
\[\pi_*L_{K(2)}(E_{\cC}^{\wedge k} \smsh X) \cong \Hom^c(\G_{\cC}^{k-1}, (E_{\cC})_*X),\]
which allows us to identify the $E_2$--term as $H^*(\G_{\cC}, (E_{\cC})_*X)$. Now, using the fact that $(E_{\cC})_*X \cong (E_{\cC})_*E_{\cC}^{hG}$ as Morava modules, (\ref{eqn:ECmod}) and Shapiro's lemma imply that  $H^*(\G_{\cC}, (E_{\cC})_*X) \cong H^*(G, (E_{\cC})_*)$. 
\end{proof}

\begin{lem}\label{lem:d3Del}
Let $X$ be a $K(2)$--local spectrum such that $(E_{\cC})_*X \cong (E_{\cC})_*E_{\cC}^{hG_{24}}$ as Morava modules. Then, the $K(2)$--local, $E_{\cC}$--Adams spectral sequence computing $\pi_{*}(X\smsh V(0))$ has $E_2$--term 
isomorphic to
$H^*(G_{24}, (E_{\cC})_*V(0))$. Further, in this spectral sequence, the elements $\Delta^k$ and $v_1\Delta^k$ are $d_3$--cycles for all $k$.
\end{lem}
\begin{proof}
The identification of the $E_2$--term follows from \fullref{lem:sse2} and the five lemma.
There are no $d_2$--differentials, so all elements survive to the $E_3$--term. 
Let $\epsilon = 0,1$. It follows from \cite[Theorem 4.2.2]{Paper1}, that $d_3(v_1^{\epsilon}\Delta^k) = v_1^{10+\epsilon}\eta^3 p(j)\Delta^{k-1}$ for $p(j) \in  \mathbb{F}_4[\![j]\!]$.
Suppose that $p(j)$ is not zero. Then $p(j) = j^rp_0(j)$ for $r\geq 0$ and $p_0(j) \in \F_4[\![j]\!]$ such that $p_0(j) \equiv \ell$ modulo $(j)$ for some $\ell \in \F_4^{\times}$. 
Using the fact that the differentials are $\eta$ and $v_1^4$--linear (since $X \smsh V(0)$ has a $v_1^4$--self map), \fullref{rem:unid3} and the identity $j = v_1^{12}\Delta^{-1}$, we have
\begin{align*}
0=d_3(v_1^{10+\epsilon}\eta^3 p(j)\Delta^{k-1} ) &=  v_1^{12r+8}\eta^3 d_3(v_1^{2+\epsilon }p_0(j)\Delta^{k-r-1}) \\
&=   v_1^{12r+8+\epsilon}\eta^6p_0(j)\Delta^{k-r-1}  + v_1^{12r+10}\eta^3d_3(v_1^{\epsilon}p_0(j)\Delta^{k-r-1}) .
\end{align*}
Again, by  \cite[Theorem 4.2.2]{Paper1}, $H^3(G_{24}, (E_{\cC})_{t}V(0))$ is $\F_4[v_1,\eta]$--torsion free in degrees $t \equiv 6+2\epsilon$ modulo $(24)$, so we can conclude that
\[ \eta^3 p_0(j)\Delta^{k-r-1}  = v_1^{2-\epsilon} d_3(v_1^{\epsilon}p_0(j)\Delta^{k-r-1}).\]
Since $\epsilon=0$ or $1$, the right hand side is divisible by $v_1$, while the left hand side is not, a contradiction. Therefore, we must have $p(j)=0$.
\end{proof}

In the next two propositions, we let
\[R(-) = \W(\!( - )\!)[\beta^{\pm 1}, \eta]/(2\eta, \eta^3).\]

\begin{prop}\label{piG24}
Let $X$ be as in \fullref{lem:d3Del}. The $E(1)$--localization of the $K(2)$--local, $E_{\cC}$--Adams spectral sequence 
\[  E_2^{s,t} = v_1^{-1}H^s(\G_{\cC}, (E_{\cC})_t(X\smsh V(0))) \Longrightarrow \pi_{t-s}L_1(X\smsh V(0))\]
satisfies
\[ E_{\infty}^{s,t}  \cong R(j)\{x, v_1x\}/( 2 \cdot x, 2v_1x)\]
for $x$ in $(0,0)$ and $v_1x\in (0,2)$.
Further, $\pi_{8t}L_1(X\smsh V(0))  \cong \F_4(\!(j)\!)\{\beta^t\}$ and the edge homomorphisms
\[  h: \pi_{8t}L_1(X\smsh V(0)) \to v_1^{-1}H^0(G_{24}, (E_{\cC})_{8t}V(0)) \]
are isomorphisms.
\end{prop}

\begin{proof}
By \fullref{lem:d3Del} and naturality, $E_2^{s,t}$ is isomorphic to $v_1^{-1}H^s(G_{24}, (E_{\cC})_tV(0))$ and
$j^k = v_1^{12k}\Delta^{-k}$ and $v_1j^k$ are $d_3$--cycles. By \fullref{rem:unid3}, there are differentials $d_3(v_1^2j^k) = \eta ^3j^k$ and $d_3(v_1^3j^k) = v_1\eta ^3j^k$. This, together with the fact that the differentials are $v_1^4$--linear, determines all $d_3$--differentials. The $E_4$--term has a horizontal vanishing line at $s=3$. Therefore, there cannot be any higher differentials. Letting $x$ be the element detected by $1 \in  H^0(G_{24}, (E_{\cC})_0V(0))$, $v_1x$ be the element detected by $v_1 \in  H^0(G_{24}, (E_{\cC})_2V(0))$ and $\beta^{t}$ be the element detected by $v_1^{4t}$, we obtain the desired description of the $E_{\infty}$--term. For degree reasons, $\pi_{8t}L_1(X\smsh V(0))  \cong \F_4(\!(j)\!)\{\beta^t\}$. That the edge homomorphisms are isomorphisms in degrees $8t$ follows since $v_1^{-1}H^0(G_{24}, (E_{\cC})_{8t}V(0)) \cong  \F_4(\!(j)\!)\{v_1^{4t}\}$ and $h(j^k\beta^{t}) = j^kv_1^{4t}$.
\end{proof}

\begin{rem}\label{rem:2v1eta2}
When $X = V(0)$, the class $x$ can be described as the composite $S^0 \to L_1E_{\cC}^{hG_{24}}  \xra{1 \smsh e_0}  L_1(E_{\cC}^{hG_{24}} \smsh V(0))$, where the first map is the unit and $e_0$ is the inclusion of the bottom cell. In $\pi_*V(0)_{(2)}$, there is a relation $2v_1e_0 = \eta^2 e_0$, for $v_1e_0$ detected by $v_1 \in BP_2V(0)$ in the Adams-Novikov spectral sequence. This then implies that $2v_1x =\eta^2 x$ in $\pi_*L_1(E_{\cC}^{hG_{24}}\smsh V(0))$, so that
\[\pi_*L_1(E_{\cC}^{hG_{24}}\smsh V(0)) \cong R(j)\{x, v_1x\}/( 2 \cdot x, 2v_1x-\eta^2x).\]

With some work, one can show that the relation $2v_1x=\eta^2x$ holds for arbitrary $X$ satisfying the condition of \fullref{lem:d3Del}. However, this fact is not needed here.
\end{rem}

\begin{figure}[H]
 \captionsetup{width=\textwidth}
\includegraphics[width=\textwidth]{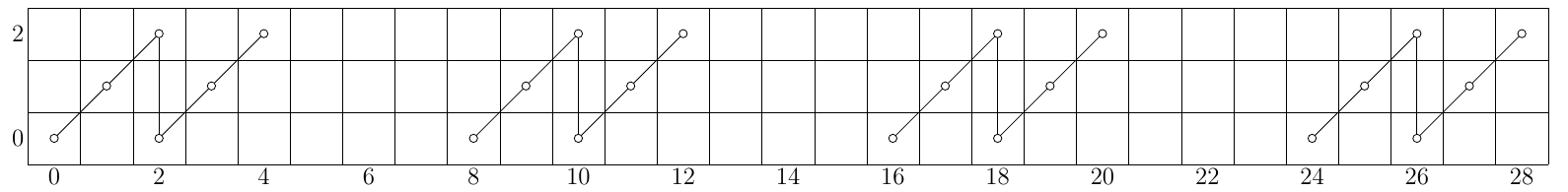}
\caption{This picture is both an illustration of the homotopy groups $\pi_*L_1(E_{\cC}^{hG_{24}}\smsh V(0))$ and of the homotopy groups $\pi_*L_1(E_{\cC}^{hC_{6}}\smsh V(0)) $. For $\pi_*L_1(E_{\cC}^{hG_{24}}\smsh V(0))$, a $\circ$ denotes a copy of $\F_4(\!(j)\!)$ and for $\pi_*L_1(E_{\cC}^{hC_{6}}\smsh V(0))$, it denotes a copy of $\F_4(\!(j_0)\!)$.}
\label{fig:C6G24top}
\end{figure}

\begin{prop}\label{piC6}
There is an isomorphism
\[\pi_*L_1(E_{\cC}^{hC_{6}}\smsh V(0)) \cong R(j_0)\{y, v_1y\}/( 2 \cdot y, 2v_1y-\eta^2y).\]
for $y$ in $(0,0)$ and $v_1y\in (0,2)$ (see Figure \ref{fig:C6G24top}). Hence, $\pi_*L_1(E_2^{hC_6}\smsh V(0))$ is $8$--periodic with periodicity generator $\beta$. Further, the edge homomorphisms
\[  h: \pi_{8t}L_1(E_{\cC}^{hC_{6}}\smsh V(0))  \to v_1^{-1}H^0(C_{6}, (E_{\cC})_{8t}V(0)) \]
are isomorphisms.
\end{prop}
\begin{proof}
We prove that $j_0^k$ is a $d_3$--cycle for all integers $k$. Then, an argument similar to that of \fullref{piG24} finishes the computation of the $E_{\infty}$--term, where we let $y$ be the element detected by $1 \in H^0(C_{6}, (E_{\cC})_0V(0))$ and $v_1y$ be the element detected by $v_1y \in H^0(C_{6}, (E_{\cC})_2V(0))$. The extension is obtained as in \fullref{rem:2v1eta2}.

The spectral sequence $H^*(C_6, (E_{\cC})_*) \Rightarrow  \pi_{*}E_{\cC}^{hC_6}$ is multiplicative, hence, in this spectral sequence, all elements of the form $a^2$ are $d_3$ cycles. Note that $j_0$ lifts to an invariant in $H^0(C_6, (E_{\cC})_0)$. This implies that $d_3(j_0^{2r})=0$ and $d_3(j_0^{2r+1}) = j_0^{2r}d_3(j_0)$. Hence, it suffices to prove that $j_0$ is a $d_3$--cycle.
The restriction induced by the inclusion of $C_6$ in $G_{24}$, maps $j$ to $j_0^4(1+j_0)^{-3}$. By naturality, the element $d_3(j_0^{4}(1+j_0)^{-3})=0$.
However,
\begin{align*}
d_3(j_0^{4}(1+j_0)^{-3}) = j_0^{4}(1+j_0)^{-4} d_3(1+j_0) = j_0^{4}(1+j_0)^{-4} d_3(j_0) ,
\end{align*}
which implies that $d_3(j_0)=0$.
\end{proof}

\section{Some Elements in $\pi_*L_1L_{K(2)}V(0)$}\label{sec:main}
We now turn to examining the spectral sequence
\[\xymatrix{E_1^{p,q}(L_1SS_1) = \pi_{q}L_1(\sE_p\smsh V(0))  \ar@{=>}[r]^-{L_1SS_1} & \pi_{q-p}L_1(E_{\cC}^{h\mathbb{S}_{\cC}^1}\smsh V(0))}. \]
The idea is to use knowledge of the differentials in the spectral sequence
 \begin{equation*}
 \xymatrix{E_1^{p,q,t}(L_1SS_2) = v_1^{-1}\Ext^q_{\Z_2[\![\mathbb{S}_{\cC}^1]\!]}(\sC_p,( E_{\cC})_tV(0) ) \ar@{=>}[r]^-{L_1SS_2} & v_1^{-1}H^{p+q}(\mathbb{S}_{\cC}^1.(E_{\cC})_tV(0))}\end{equation*}
to deduce information about the differentials of $L_1SS_1$.
\begin{lem}\label{lem:e23}
In the spectral sequence $L_1SS_1$, $E_2^{3,8t} \cong \F_4\{\beta^t\}$.
\end{lem}
\begin{proof}
From \fullref{sec:e1page}, we have that
\[ E_1^{p,8t} \cong \begin{cases}
\F_4(\!(j )\!)\{\beta^t\} & p=0,3 \\ 
\F_4(\!(j_0 )\!)\{\beta^t\} & p=1,2.
\end{cases}\]
From \fullref{rem:edgehoml1} and the fact that the edge homomorphisms are isomorphisms in these degrees, we obtain a commutative diagram
\[\xymatrix{E_1^{0,8t}(L_1SS_1) \ar[r]^-{d_1}  \ar[d]^-{\cong}  & E_1^{1,8t}(L_1SS_1)  \ar[r]^-{d_1}   \ar[d]^-{\cong}  & E_1^{2,8t}(L_1SS_1) \ar[r]^-{d_1}  \ar[d]^-{\cong}  & E_1^{3,8t}(L_1SS_1) \ar[d]^-{\cong} \\
 E_1^{0,0, 8t}(L_1SS_2)  \ar[r]^-{d_1}  &  E_1^{1,0,8t}(L_1SS_2)  \ar[r]^-{d_1} &  E_1^{2,0,8t}(L_1SS_2)  \ar[r]^-{d_1} &  E_1^{3,0,8t}(L_1SS_2),   }\]
 where $\beta^{4t}$ maps to $v_1^{4t}$.
Theorem 1.2.1 and Corollary 1.2.3 of \cite{Paper2} give a computation of the spectral sequence $L_1SS_2$. In particular, it follows immediately from these results that
\[E_2^{3,0,8t}(L_1SS_2) \cong \F_4(\!(j )\!)\{v_1^{4t}\}/(j) \cong \F_4\{v_1^{4t}\} . \]
The claim follows.
\end{proof}

\begin{prop}\label{pi5S21}
If $k \equiv -3$ modulo $8$, then $\pi_{k}L_1(E_{\cC}^{h\mathbb{S}_{\cC}^1}\smsh V(0)) \cong \F_4$.
\end{prop}
\begin{proof}
We use the spectral sequence $E_r^{p,q} = E_r^{p,q} (L_1SS_1)$. From \fullref{piG24} applied to $X = \sE_0$ and $X = \sE_3$ and from \fullref{piC6}, it follows that for $r=1$, $2$ or $3$ and for any $p$,
\[E_1^{p, 8t-r} = \pi_{8t-r}L_1(\sE_p \smsh V(0)) =0.\] 
By \fullref{lem:e23}, $E_2^{3,8t} \cong \F_4\{\beta^{8t}\}$, which proves the claim.
\end{proof}

\begin{prop}\label{pi5G21}
If $k \equiv -3$ modulo $8$, then $\pi_{k}L_1(E_{\cC}^{h\G_{\cC}^1}\smsh V(0)) \cong \F_2$. 
\end{prop}
\begin{proof}
It follows from \fullref{rem:Wlinl1} that $v_1^{-1}E_{\infty}^{*,*}(\Sn_{\cC}^1,V(0)) \cong   \W \otimes_{\Z_2} v_1^{-1}E_{\infty}^{*,*}(\G_{\cC}^1,V(0))$. Since $\pi_{k}L_1(E_{\cC}^{h\mathbb{S}_{\cC}^1}\smsh V(0)) \cong \F_4$, there exists a unique $s_0\geq 0$ such that $E_{\infty}^{s_0,k+s_0}(\Sn_{\cC}^1,V(0))$ is non-zero, and $E_{\infty}^{s_0,k+s_0}(\Sn_{\cC}^1,V(0))\cong \F_4$. Therefore, $E_{\infty}^{s,k+s}(\G_{\cC}^1,V(0)) =0$ if $s\neq s_0$ and $E_{\infty}^{s_0,k+s_0}(\G_{\cC}^1,V(0)) \cong\F_2$.
\end{proof}

\begin{defn}\label{defrho}
We define the class $x \in \pi_{-3}L_1(E_{\cC}^{h\G_{\cC}^1}\smsh V(0))$ to be the non zero element.
\end{defn}

Recall that
\[\G_{\cC}  \cong \G_{\cC}^1 \rtimes \Z_2.\]
Let $\pi$ be a topological generator of the subgroup $\Z_2$ in $\G_{\cC}$. There is a cofiber sequence
\begin{equation}\label{fibseq}
L_{K(2)}S \ra E_{\cC}^{h\G_{\cC}^1}  \xra{\pi-1} E_{\cC}^{h\G_{\cC}^1}.
\end{equation}
We can now prove our main result.
\begin{proof}[Proof of \fullref{main}.]
Since $L_{K(2)}S\smsh V(0)\simeq L_{K(2)}V(0)$ and localization preserves exact triangles,
 the fiber sequence (\ref{fibseq}) gives rise to a fiber sequence
\begin{equation}
L_1L_{K(2)}V(0) \ra L_1(E_{\cC}^{h\G_{\cC}^1}\smsh V(0)) \xra{\pi-1} L_1(E_{\cC}^{h\G_{\cC}^1}\smsh V(0)).
\end{equation}
Since $\pi$ acts by automorphisms and the only automorphism of $\F_2$ is the identity, the map $\pi-1$ acts trivially on $\pi_{8t-3}L_1(E_{\cC}^{h\G_{\cC}^1}\smsh V(0)) $. Therefore, in the long exact sequence on homotopy groups, the class $\beta^t x$ is in the kernel of $\pi-1$, and the image of $\beta^t x$ under the map $L_1(E_{\cC}^{h\G_{\cC}^1}\smsh V(0)) \ra \Sigma L_1L_{K(2)}V(0) $ is non-zero. We denote it by $\zeta_2 \beta^t x$. 
\end{proof}


%
%
%
%

\end{document}